\documentclass[10pt]{article}
\usepackage{arXiv}
\usepackage{thmtools}
\usepackage{thm-restate}

\usepackage{todonotes}

\usepackage{fancyhdr}
\title{
    \textbf{
    Accelerating Multivariate Newton Interpolation in Downward Closed Polynomial Spaces}
}

\author{ 
    Phil-Alexander Hofmann\thanks{
        Center for Advanced Systems Understanding, 
        Helmholtz-Zentrum Dresden-Rossendorf e.V.,
        Untermarkt 20, 02826 Görlitz, 
        Germany, 
        \texttt{p.hofmann@hzdr.de} }, \
    Michael Hecht\thanks{ 
        Mathematical Institute, 
        University of Wrocław, 
        pl. Grunwaldzki 2/4, 
        50-384 Wrocław, 
        Poland, 
        \texttt{michael.hecht@math.uni.wroc.pl} 
    } 
}

\begin{document}

\maketitle
\begin{abstract}
    We introduce the fast Newton transform (FNT), a multivariate Newton interpolation algorithm for downward closed polynomial spaces in quasi-tensorial grids.
    The FNT computes the Newton coefficients directly, without relying on embeddings into enclosing tensor-product spaces. 
    For a downward closed index set $A \subset \mathbb N_0^m$, the FNT achieves a time complexity of $\mathcal O(m \overline n |A|)$, where $\overline n$ is the mean of the coordinate-wise maximal polynomial degrees $n_1, \ldots, n_m$ across the $m$ spatial dimensions. 
    In the univariate case, the FNT renders the classic Newton divided difference scheme (DDS). 
    In the multivariate case, however, it improves on the quadratic complexity $\mathcal O(|A|^2)$ of the DDS and on the cost $\mathcal{O}(m \overline n (n_1+1) \cdots (n_m+1))$ of the tensorial interpolation.
    
    For sufficiently regular functions, Newton interpolation in Euclidean-degree downward closed polynomial spaces is known to deliver approximation rates equal to those of the tensor product interpolation. 
    Thus, the acceleration power of the FNT comes from requiring substantially fewer degrees of freedom while reaching the same approximation quality as the tensorial interpolation. 
    The inverse transformation has the same time complexity and enables fast evaluation and differentiation of Newton interpolants in quasi-tensorial grids.
\end{abstract}

\section{Introduction}

Across disciplines, the computational task of expanding multivariate functions into closed-form expressions is omnipresent. This includes solving partial differential equations \cite{Gottlieb1987, Boyd2000, Press2007}, uncertainty quantification \cite{Sudret2008, LeMaitre2010, Adelmann2019, Castellon2023, Loukrezis2025}, optimization \cite{Fischer2011, Lasserre2015, Ayvaz2022, Schreiber2023}, neural networks \cite{Suarez2024, Morala2025}, computer-aided design \cite{Speleers2011, Veettil2023}, and inverse problems \cite{Fatih2013, Arridge2019}.

Interpolation of multivariate non-periodic functions in downward closed polynomial spaces was suggested and studied in recent years  \cite{Chkifa2013, Chkifa2014, Chkifa2015, Cohen2015, Trefethen2017a, Trefethen2017b, Cohen2017, Cohen2018, Bos2018, Hecht2026}. In particular, for regular functions, suitable choices of the downward closed polynomial space mitigate the curse of dimensionality effectively. We briefly outline this perspective.

\subsection{Related work}
The classic theorems of Weierstrass, Bernstein, and Stone prove the capabilities of polynomials to approximate continuous functions to be universal~\cite{Weierstrass1885, Bernstein1912, Stone1959}. In high dimensions, however, the curse of dimensionality~\cite{Bellman1957, Bellman1961} demands sparse expansions, in order to reduce data acquisition and computational costs, see, e.g., \cite{Babuska1997, Xiu2002, Dick2013, Hu2024}.

Established approaches for multivariate interpolation and other related numerical tasks include sparse grid methods \cite{Smolyak1963, Gerstner1998, Bungartz2004}, which exploit tensor-product structure and mixed smoothness, and radial basis function methods \cite{Broomhead1988, Park1991}, which use superpositions of radially symmetric kernels. Here, we follow a polynomial approach based on downward closed index sets. In this setting, analytic continuation to an open Bernstein poly-ellipse provides a sufficient condition for geometric approximation rates \cite{Trefethen2017a, Bos2018, Hecht2026}.

The explicit term \emph{downward closed polynomial spaces} was formalized by \emph{Cohen} and \emph{Migliorati}~ \cite{Cohen2017}. They form linear polynomial subspaces that are induced by \emph{downward closed index sets} and closed under differentiation. 
The alternative terminology is \emph{lower sets}, a notion used in lattice and order theory~\cite{Priestley2002}.

\emph{Chkifa et al.} emphasized the numerical practicality of using downward closed index sets for anisotropic polynomial approximations, demonstrating that algebraic convergence rates of tensorized Legendre expansions in infinite-dimensional settings are preserved for polynomial interpolation \cite{Chkifa2013, Chkifa2014, Chkifa2015}.

When it comes to the isotropic case, a key observation was made by \emph{Trefethen}. That is for functions extending analytically in an unbounded set, constructed via the \emph{Newton ellipse}~\cite{Kazarinoff1991}, the Euclidean degree polynomial spaces \(\Pi_{m,n,2}\), deliver the same approximation power as the maximum degree spaces \(\Pi_{m,n,\infty}\), while the total degree \(\Pi_{m,n,1}\) fail to resolve the function along the diagonal of \(A_{m,n,1}\) by a factor of \(\sqrt{m}\)~\cite{Trefethen2017a, Trefethen2017b}. 

Meanwhile, \emph{Bos} and \emph{Levenberg} have further deepened the approximation theory,  establishing a Bernstein-Walsh Theorem for downward closed polynomial spaces generated by convex bodies \cite{Bos2018}. 
\emph{Hecht et al.} \cite{Hecht2026} demonstrated that these analytic functions can be expanded by multivariate \emph{Newton} interpolation in Leja ordered \emph{Chebyshev}-\emph{Lobatto} nodes or \emph{Leja} points \cite{Leja1957, Andrievskii2022} with the predicted optimal rates. This article contributes by accelerating the runtime of the former interpolation algorithm, as we outline next.

\subsection{Contribution}

Algorithms for multivariate polynomial interpolation have been studied extensively \cite{Muehlbach1978, Brezinski1980, Muehlbach1988, Boor1992, Sauer1995a, Sauer1995b}. In particular, the DDS has been extended to this setting \cite{Boor1995, Olver2006, Neidinger2019}. 
More recently, \emph{Chkifa et al.} and \emph{Cohen et al.} \cite{Chkifa2014, Cohen2018} employed an algorithm for interpolation in downward closed polynomial spaces $\Pi_A$, $\dim \Pi_A = |A|$, requiring super-quadratic \(\Omega(|A|^2)\) up to cubic runtime \(\mathcal{O}(|A|^3)\). \emph{Hecht et al.} \cite{Hecht2026} presented a \emph{Newton} interpolation algorithm with quadratic runtime \(\mathcal{O}(|A|^2)\).

We present the \emph{fast Newton transform} (FNT), an algorithm that computes the \emph{Newton} coefficients in \(\mathcal{O}(m\overline{n} |A|)\), see Theorem~\ref{theo:fnt}. The inverse FNT has the same time complexity and realizes evaluation in quasi-tensorial grids. Moreover, the \emph{Newton} coefficients of partial derivatives can be computed in \(\mathcal{O}(n_i |A|)\), see Theorem~\ref{theo:fastdiff}, and subsequently evaluated in the grid via the inverse FNT in \(\mathcal{O}(m\overline n|A|)\).

Specifically, we exploit a Kronecker-type product representing the \emph{Vandermonde} matrix and its inverse for general separable polynomial bases of downward closed polynomial spaces, see Proposition~\ref{prop:selected_kronecker} and Proposition~\ref{prop:selected_kronecker_inverse}.
The key observation is that, for \emph{Newton} interpolation, the triangular structure of the univariate factors is feasible for downward closed sets. 
This allows the Kronecker-type product to be decomposed into mode-wise factors, see Proposition~\ref{prop:mode_factorization}, and yields the FNT.
In general, for other bases, the identified Kronecker-type product can not be evaluated directly by successive univariate factors.
Herefore, univariate LU decompositions of the corresponding \emph{Vandermonde} matrices have to be computed a priori passing the computation again through the \emph{Newton} basis.
To the best of our knowledge, this structure has not been observed and exploited prior for interpolation in downward closed spaces.

A key contribution for practitioners is the implementation of the \emph{FNT} for a family of \(\ell^p\)-type index sets in the open-source Python package \emph{lpFun}~\cite{Hofmann2026lpFun}, which leverages \emph{Numba} for C-level performance. 

\subsection{Notation and preliminaries}

We aim to be as consistent as possible with the conventions established in the field of tensor decompositions. The notation closely follows the conventions proposed by \emph{Kolda} and \emph{Bader}~\cite{Kolda2009}, originating from \emph{Kiers}~\cite{Kiers2000}. For downward closed polynomial spaces, we follow the notation used by \emph{Hecht et al.}~\cite{Hecht2026}.

\emph{Scalars} are denoted by lowercase letters, for example \(x\). 
\emph{Vectors} are denoted by boldface lowercase letters, for example \(\bm x\). 
\emph{Matrices} are denoted by boldface capital letters, for example \(\bm X\). 
\emph{Higher order tensors} are denoted by boldface calligraphic letters, for example \(\bm{\mathcal X}\).
\emph{Transformation matrices} are denoted by capital Greek letters, for example \(\Psi\). 
\emph{Multi-indices} are denoted by boldface lowercase Greek letters, for example \(\bm \alpha\). 
\emph{Multi-index sets} are denoted by capital letters, for example \(A\).

The \(i\)th entry of a vector \(\bm x\) is denoted \(x_i\). If the vector \(\bm x\) is indexed by a multi-index set \(A\), we equivalently write its entries as \(x_{\bm\alpha}\), \(\bm\alpha\in A\).
The element in the \(i\)th row and \(j\)th column of a matrix \(\bm X\) is denoted by \(x_{i,j}\).
The element of a tensorized matrix \(\bm{\mathcal X}\), indexed by multi-indices \(\bm\alpha\) and \(\bm\beta\), is denoted by \(x_{\bm\alpha,\bm\beta}\).
The \(n\)th element in a sequence is denoted by a superscript in parentheses, for example \(\bm X^{(n)}\) denotes the \(n\)th matrix in a sequence.

For convenience, Table~\ref{tab:notation} serves as a complementary reference listing the notation and basic operations. \\

% Notation table
\begin{table}[ht]
    \centering
    \caption{Notation and basic operations.}
    \label{tab:notation}
    \renewcommand{\arraystretch}{1.15}
    \setlength{\tabcolsep}{8pt}
    \small
    \begin{tabular}{|c|l||c|l||c|l|}
    \hline
    $\bm x$ & vector & $\bm X$ & 
    matrix & $\bm{\mathcal X}$ & tensor \\
    $A$ & multi-index set & 
    $\bm\alpha, \bm\beta, \bm\gamma$ & multi-indices &
    $\Pi_A$ & polynomial space \\
    $|\cdot|$ & set cardinality &
    $\|\cdot\|_p$ & $\ell^p$-norm &
    $\|\cdot\|_{L^p}$ & \(L^p\) norm\\
    $\mathcal{O}(\cdot)$ & asymptotic upper bound & 
    $\Omega(\cdot)$ & asymptotic lower bound & 
    $C(\Omega)$ & continuous functions\\
    $\operatorname{span}$ & linear hull &
    $\bm I^{(n)}$ & $n\times n$ identity matrix &
    $S_n$ & Symmetric group \\
    $\times$ & Cartesian product &
    $\otimes$ & Kronecker product &
    $\Lambda(\cdot)$ & Lebesgue constant\\
    \hline
    \end{tabular}
\end{table}

\vspace{1em}

The subsequent Sections~\ref{sec:interpolation_on_dcps} and \ref{sec:kronecker_product_structure} develop the computational framework leading to the concrete \emph{FNT} algorithm in Section~\ref{sec:construction_of_the_fnt}.

\section{Interpolation in downward closed polynomial spaces}
\label{sec:interpolation_on_dcps}

We formalize the ingredients on which the present interpolation setup rests.

\subsection{Downward closed polynomial spaces}
\label{subsec:dcps_and_indexing}

We begin with defining the notion of downward closed polynomial spaces.

\begin{definition}[Downward closed polynomial space]
    \label{def:dcps}
    Let \(m \in \mathbb N\) and \(A \subset \mathbb N_0^m\). \(A\) is downward closed, if
    \begin{equation}
        \forall\,\bm\alpha = (\alpha_1, \ldots, \alpha_m)\in A \,: \, \{\bm\beta = (\beta_1, \ldots, \beta_m) \in \mathbb N_0^m \, : \, \beta_1 \leq \alpha_1, \ldots, \beta_m\leq \alpha_m\} \subset A.
    \end{equation}
    Then, we define the downward closed polynomial space with respect to \(A\) by
    \begin{equation}
        \Pi_A \coloneqq \operatorname{span}_\mathbb R \{\bm x^{\bm\alpha} \coloneqq x_1^{\alpha_1} \cdots  x_m^{\alpha_m} \, : \, \bm\alpha \in A\}.
    \end{equation}
    If the set \(A\) is finite, we denote the maximal polynomial degree of each spatial dimension, the average and the overall maximum respectively by
    \begin{equation}
        n_{i} \coloneqq \max_{\alpha \in A} \alpha_i, \qquad 
        \overline{n} \coloneqq \frac{n_1 + \cdots + n_m}{m}, 
        \qquad 
        n \coloneqq \max_{i=1,\ldots,m} n_i.
    \end{equation}
\end{definition}

\begin{remark}
    A natural choice for a total order on downward closed index sets is the \emph{lexicographic order}. It arranges sequences by comparing their elements from the first to the last, in other words, in dictionary order. For elements
    \(
        \bm\alpha, \bm\beta \in \mathbb N_0^m,
    \)
    we define
    \begin{equation}
        \bm\alpha \preceq_{\text{lex}} \bm\beta \iff \bm \alpha = \bm \beta \text{ or } \exists \, i \in \{1,\ldots,m\}: \alpha_i < \beta_i, \, \alpha_{i-1}=\beta_{i-1}, \ldots, \alpha_1 = \beta_1.
    \end{equation}
\end{remark}

In particular, we switch between indexing vectors by \(\bm\alpha \in A \subset \mathbb N_0^m\) and by their ordinal positions \(k = 1, \ldots, |A|\). To formalize this, we consider the rank map, which will also be used later.

\begin{definition}[Rank map]
    Let \(m \in \mathbb N\), let \(A \subset \mathbb N_0^m\) be finite, downward closed and equipped with a total order \(\preceq\). We denote the ordinal position of \(\bm\alpha \in A\) by
    \begin{equation}
        \operatorname{rank}_A(\bm\alpha) 
        \coloneqq 
        |\{\bm\beta \in A \ : \ \bm\beta \preceq \bm\alpha\}|.
    \end{equation}
    The inverse of \(\operatorname{rank}_A\) exists, allowing us to retrieve \(\bm\alpha \in A\) from its ordinal position. Hence,
    \(
        \operatorname{rank}_A^{-1}(k)
    \) 
    for \(k \in \{1, \ldots, |A|\}\), denotes the \(k\)th element of \(A\).
\end{definition}

\subsection{Quasi-tensorial grids and Vandermonde matrices}
\label{subsec:quasi_tensorial_grids_and_vandermonde_matrices}

Following \emph{Hecht et al.}~\cite{Hecht2026}, we consider the hypercube \([-1, 1]^m\) and define quasi-tensorial grids used for interpolation in downward closed polynomial spaces.

\begin{definition}[Quasi-tensorial grid]
    \label{def:quasi_tensorial_grid}
    Let \(m \in \mathbb N\), let \(A \subset \mathbb N_0^m\) be finite and downward closed. Further, let arbitrary sets of pairwise distinct nodes 
    \(
        P_i \coloneqq \left\{p_{0,i}, p_{1,i}, \ldots, p_{n_i,i}\right\} \subset [-1, 1].
    \)
    We consider the quasi-tensorial grid \(P_A\) associated with \(A\), defined as
    \begin{equation}
        \label{eq:quasi_tensorial_grid}
        P_A \coloneqq \left\{\bm p_{\bm\alpha} 
        \coloneqq 
        (p_{\alpha_1,1},\ldots,p_{\alpha_m,m}) \, : \, \bm\alpha \in A\right\}.
    \end{equation}
\end{definition}

Associated with this quasi-tensorial grid \(P_A\), the following \emph{Vandermonde} matrix represents the interpolation problem in \(\Pi_A\).

\begin{definition}[Vandermonde matrix]
    \label{def:vandermonde_matrix}
    Let \(m \in \mathbb N\), let \(A \subset \mathbb N_0^m\) be finite and downward closed. Further, let \(P_A\) be a quasi-tensorial grid and \(Q_A \coloneqq \{Q_{\bm\beta} \, : \, \bm\beta \in A\}\) be a basis of \(\Pi_A\). 
    We consider the \emph{Vandermonde} matrix
    \begin{equation}
        \bm{\mathcal V} \coloneqq (Q_{\bm\beta}(\bm p_{\bm\alpha}))_{\bm\alpha, \bm\beta \in A} \in \mathbb R^{|A|\times|A|}.
    \end{equation}
    The polynomial basis \(Q_A\) is called \emph{separable}, if there exist polynomials \(Q_{j,i} \in \mathrm{span}_{\mathbb R} \left\{x^0, x^1, \ldots, x^{n_i}\right\}\), for all \(j=0,\ldots,n_i\), \(i=1,\ldots,m\), such that 
    \begin{equation}
        \label{eq:separable_polynomial_basis}
         Q_{\bm\beta}(\bm x) 
         = 
         Q_{\beta_1,1}(x_1) \cdots Q_{\beta_m,m}(x_m),
    \end{equation}
    for all \(\bm\beta \in A\). Then, if \(Q_A\) is separable, we denote
    \begin{equation}
        \label{eq:univariate_vandermonde_matrix}
        \bm V^{(i)} \coloneqq \left(Q_{k,i}(p_{j,i})\right)_ {j,k=0,\ldots,n_i}.
    \end{equation}
\end{definition}

The \emph{Vandermonde} matrices \(\bm V^{(i)}\) are invertible if and only if their underlying nodes are distinct, leading to unisolvent grids \(P_A\), see \cite{Hecht2026}.

\subsection{Polynomial interpolation schemes}
\label{subsec:separable_interpolation_schemes}

We now formalize the interpolation problem in \(\Pi_A\). 
Our terminology distinguishes between general separable interpolation and separable \emph{Newton} interpolation.
In a \emph{separable interpolation scheme}, the multivariate structure is induced by univariate bases and nodes. 
In a \emph{separable Newton interpolation scheme}, the basis is, in addition, specifically chosen as the multivariate \emph{Newton} basis.
This additional requirement renders the associated univariate \emph{Vandermonde} matrices to be lower triangular, which is the key ingredient for interpolation in downward closed polynomial spaces.

\begin{definition}[Multivariate Newton polynomials]
    \label{def:multivariate_newton_polynomials}
    Let \(m \in \mathbb N\), let \(A \subset \mathbb N_0^m\) be finite, downward closed and lexicographically ordered. We define the \emph{m-variate Newton polynomials} associated with \(A\) to be the collection \(Q_A\) consisting of all
    \begin{equation}
        Q_{\bm\beta}(\bm x) 
        = 
        \prod_{i=1}^m\prod_{k=0}^{\beta_i-1}(x_i-p_{k,i}), 
        \quad 
        \bm \beta \in A,
    \end{equation}
    with the convention that \(\prod_{k=0}^{-1}(\ldots) \coloneqq 1\).
\end{definition}

This basis choice specifies the \emph{Newton} case within the interpolation schemes below.

% The corresponding univariate matrices \(\bm V^{(i)}\) as defined in \eqref{eq:univariate_vandermonde_matrix} can be computed dimension-wise in \(\mathcal O(n_i^2)\) operations via the Horner scheme, while their inverses \(\big(\bm V^{(i)}\big)^{-1}\) can be obtained by divided differences in \(\mathcal O(n_i^3)\).

\begin{definition}[Polynomial interpolation schemes]
    \label{def:interpolation_schemes}
    Let \(m \in \mathbb N\), let \(A \subset \mathbb N_0^m\) be finite and downward closed. Further, let \(P_A\) be a quasi-tensorial grid and \(Q_A \coloneqq \{Q_{\bm\beta} \ : \ \bm\beta \in A\}\) be a basis of \(\Pi_A\). For \(f \in C([-1,1]^m)\), polynomial interpolation with respect to \((Q_A,P_A)\) delivers the unique polynomial \(\mathcal Q_f \in \Pi_A\) coinciding with \(f\) in \(P_A\)
    \begin{equation}
        \label{eq:polynomial_interpolation}
        \mathcal Q_f(\bm p_{\bm \alpha})
        = 
        f(\bm p_{\bm\alpha}),
        \qquad 
        \text{for all }\bm \alpha \in A.
    \end{equation}
    Equivalently, the coefficients \(\bm c \in \mathbb R^{|A|}\) of \(\mathcal Q_f\) are given as the  unique solution of the linear system  
    \begin{equation}
        \bm{\mathcal V} \bm c = \bm f,
        \qquad 
        \bm f \coloneqq  (f(\bm p_{\bm\alpha}))_{\bm\alpha \in A} \in \mathbb R^{|A|}
    \end{equation}
    We call the pair \((Q_A, P_A)\) together with \eqref{eq:polynomial_interpolation}
    \begin{enumerate}[label=(\roman*)]
        \item a \emph{separable interpolation scheme}, if \(Q_A\) is separable;
        \item a \emph{separable Newton interpolation scheme}, if \(Q_A\) is the multivariate \emph{Newton} basis from Definition~\ref{def:multivariate_newton_polynomials}.
    \end{enumerate}
    We identify the \emph{two main computational tasks} to be:
    \begin{enumerate}[label=\textbf{C\arabic*)}]
        \item \label{C1} Given the function values \(\bm f\), compute the coefficients \(\bm c\) of the interpolant  \(\mathcal Q_f\).
        \item \label{C2} Given the coefficients \(\bm c\), evaluate the polynomial \(\mathcal Q_f\) in all grid points \(P_A\).
    \end{enumerate}
\end{definition}

Following this, we state the main complexity results for these two tasks.

\subsection{Main algorithmic results}
\label{subsec:main_results}

The \emph{fast Newton transform} is an algorithm that generalizes fast Kronecker-product vector multiplication to downward-closed index sets by exploiting the hierarchical triangular structure of the \emph{Vandermonde} matrix and the quasi-tensorial structure of the interpolation grid, thereby avoiding embeddings into enclosing tensor-product spaces that would increase the computational cost.

\begin{restatable}[The fast Newton transform]{theorem}{fastnewtontransform}
    \label{theo:fnt}
    Let \((Q_A, P_A)\) be a separable \emph{Newton} interpolation scheme. Then, given that \(\bm{\mathcal V} \bm c = \bm f\), the \emph{Newton} transform \(\bm f \longmapsto \bm c\) corresponding to the computation of \emph{Newton} coefficients in \ref{C1}, and its inverse \(\bm c \longmapsto \bm f\), corresponding to the evaluation task in \ref{C2}, can each be computed in 
    \begin{equation}
        \mathcal{O}(m \overline n |A|),
    \end{equation}
    where \(\overline n = (n_1+ \cdots + n_m)/m\) denotes the average of the maximal degrees of \(A\).
\end{restatable}
\begin{proof}
    The proof is given in Section~\ref{sec:construction_of_the_fnt}.
\end{proof}

We deduce that Theorem~\ref{theo:fnt} enables fast differentiation in downward closed polynomial spaces.

\begin{restatable}{theorem}{fastdiff}
    \label{theo:fastdiff}
    Let the assumptions of Theorem~\ref{theo:fnt} be fulfilled. Given the coefficients of a polynomial \(q \in \Pi_A\), the coefficients of its partial derivative \(\partial_{x_i} q\) can be computed in 
    \begin{equation}
        \mathcal O(n_i |A|).
    \end{equation}
\end{restatable}
\begin{proof}
    The proof is also given in Section~\ref{sec:construction_of_the_fnt}.
\end{proof}

Both Theorem~\ref{theo:fnt} and Theorem~\ref{theo:fastdiff} also extend to separable interpolation schemes, which is a direct consequence of Theorem~\ref{theo:selected_kronecker_lu}. Next, we present our implementation choices that serve as the baseline for the numerical experiments in Section~\ref{sec:numerical_experiments}.

\subsection{Implementation choices}
\label{sec:implementation_choices}

We specify the implementation choices for the interpolation scheme \((Q_A, P_A)\) implemented in \emph{lpFun}~\cite{Hofmann2026lpFun}. These choices consist of the downward closed index set \(A\), the basis \(Q_A\) of the downward closed polynomial space \(\Pi_A\) and the quasi-tensorial grid \(P_A\).

To begin with, the index sets used in \emph{lpFun} are chosen from a family of \(\ell^p\)-type downward closed index sets.

\paragraph{Choosing the index sets \(A\).} Let \(m, n \in \mathbb N\), \(p \in (0, \infty]\), and let the \(\ell^p\)-type index set
\begin{equation}
    A_{m,n,p} 
    \coloneqq 
    \{
        \bm\alpha \in \mathbb N_0^m \ : \ \|\bm\alpha\|_{p} \leq n 
    \}.
\end{equation}
This induces the downward closed polynomial space 
\begin{equation}
    \Pi_{m,n,p} \coloneqq \operatorname{span}_{\mathbb R} \{\bm x^{\bm \alpha} \ : \bm \alpha \in A_{m,n,p}\}.
\end{equation}
The parameter \(p\) controls the sparsity of the resulting index set, ranging from sparse sets for \(0 < p < 1\), through total-degree type spaces for \(p = 1\) and Euclidean-degree type spaces for \(p=2\), to maximal-degree type spaces for \(p = \infty\).

This choice is motivated by an approximation result of \emph{Bos} and \emph{Levenberg}~\cite{Bos2018}, which provides Bernstein-Walsh type rates for polynomial spaces associated with convex bodies, i.e. for \(1 \le p \le \infty\), and by the practical implementation in \emph{Hecht et al.}~\cite{Hecht2026}. Let \(f \in C([-1,1]^m)\) admit a holomorphic extension to the open \(\ell^q\) poly-ellipse
\begin{equation}
    \Omega_\rho^q \coloneqq \left\{ \bm z \in \mathbb C^m \, : \, \|(\log|z_j + \sqrt{z_j^2 - 1}|)_{j=1,...,m}\|_q < \log \rho\,\right\}, \quad \rho > 1,
\end{equation}
where \(q\) is the conjugate of \(p\) defined by \(1/p+1/q=1\), using the conventions \(q=\infty\) for \(p=1\) and \(q=1\) for \(p=\infty\). Then, for every \(\varepsilon \in (0,\rho)\), a geometric approximation rate is guaranteed
\begin{equation}
    \inf_{r \in \Pi_{m,n,p}} \|f - r\|_{L^\infty(\Omega)} 
    \in 
    \mathcal O((\rho - \varepsilon)^{-n}).
\end{equation}

\paragraph{Choosing the bases \(Q_A\).}

The basis \(Q_A=\{Q_{\bm\alpha}:\bm\alpha\in A\}\) can be specified by choosing a univariate polynomial family and then forming the corresponding separable multivariate basis. Besides the \emph{Newton} basis that was introduced in Definition~\ref{def:multivariate_newton_polynomials}, we implemented \emph{Chebyshev} and \emph{Legendre} bases. In one dimension, these are given by
\begin{equation}
    T_k(x) \coloneqq \cos (k \arccos x),
    \qquad
    L_k(x) \coloneqq \frac{1}{2^k k!} \frac{d^k}{dx^k}(x^2-1)^k,
    \qquad 
    k \in \mathbb N_0,
\end{equation}
with \(T_k\) denoting the \emph{Chebyshev} polynomial of the first kind and \(L_k\) denoting the Legendre polynomial, each of degree \(k\). Hence, the corresponding separable multivariate bases are
\begin{equation}
    T_{\bm\alpha}(\bm x)
    \coloneqq
    \prod_{i=1}^m T_{\alpha_i}(x_i),
    \qquad
    L_{\bm\alpha}(\bm x)
    \coloneqq
    \prod_{i=1}^m L_{\alpha_i}(x_i),
    \qquad
    \bm\alpha\in A.
\end{equation}
For the \emph{Newton} basis, the univariate \emph{Vandermonde} matrices are lower triangular, which allows the \emph{FNT} to be applied directly.
For the \emph{Chebyshev} and Legendre bases,  we use LU-decompositions of the univariate \emph{Vandermonde} matrices.
Without such a factorization, the full \emph{Vandermonde} matrices do not exhibit the structure which is required to apply the \emph{FNT}.

\paragraph{Choosing the grids \(P_A\).}
For generating the quasi-tensorial grids \(P_A\), we follow the construction of \emph{Hecht et al.}~\cite{Hecht2026}, which is based on the \emph{Leja ordering}~\cite{Leja1957}. This ordering is used because it provides a numerically stable ordering of interpolation nodes for high degree polynomial interpolation in \emph{Newton} form~\cite{TalEzer1991}.
Specifically, given pairwise distinct nodes \(x_1, \ldots, x_n\) with \(1 \in \{x_1,\ldots,x_n\}\), a \emph{Leja ordering} is obtained by iteratively choosing
\begin{equation}    
    y_j \in \arg\max_{y\in\{x_1,\ldots,x_n\}\setminus\{y_1,\ldots,y_{j-1}\}} \prod_{\ell=1}^{j-1}|y-y_\ell|,
    \qquad
    y_1 \coloneqq 1.
\end{equation}
Each new node is chosen to maximize its multiplicative distance from the previously selected nodes. 
Equivalently, this amounts to a greedy maximization of the \emph{Vandermonde} determinant.
In particular, if the one-dimensional nodes \(P_i\) fulfill
\begin{equation}
    \Lambda(\{p_{0,i}, p_{1,i}, \ldots, p_{k,i}\}) \le C (k+1)^\theta,
    \qquad
    i = 1, \ldots, m,
\end{equation}
for a constant \(C>0\) and \(\theta \ge 1\), the Lebesgue constant of the quasi-tensorial grid \(P_A\) satisfies
\begin{equation}
    \Lambda(P_A) \le C' |A|^{\theta+1},
\end{equation}
with a constant \(C' > 0\) independent of \(|A|\), as shown by \emph{Chkifa et al.}~\cite{Chkifa2014}.

Building upon this, we consider two choices for the one-dimensional nodes \(P_i\). 
The first choice consists of Leja sequences in \([-1,1]\), generated by applying the Leja selection procedure to a sufficiently fine discretization. 
For Leja points in \([-1,1]\), the one-dimensional Lebesgue constants satisfy the upper bound~\cite{Andrievskii2022}
\begin{equation}
    \Lambda(\{p_{0,i}, p_{1,i}, \ldots, p_{k,i}\}) \in \mathcal O ((k+1)^{13/4}).
\end{equation}
Consequently, the Lebesgue constant of \(P_A\) grows only polynomially
\begin{equation}
    \Lambda(P_A) \in \mathcal O (|A|^{17/4}).
\end{equation}
As a second choice, we use \emph{Chebyshev}-\emph{Lobatto} nodes
\begin{equation} 
    P_i = \left\{ \cos\left(k \pi / n_i\right) \ : \ k=0, \ldots, n_i \right\}, 
\end{equation} 
equipped with a \emph{Leja ordering}. Since the one-dimensional \emph{Chebyshev}-\emph{Lobatto} nodes satisfy\linebreak
\(
    \Lambda(P_i) \in \mathcal O (\log n_i),
\)
the Lebesgue constant of the first \(k+1\) nodes in Leja order are bounded by
\begin{equation}
    \Lambda(\{p_{0,i}, p_{1,i}, \ldots, p_{k,i}\})
    \in
    \mathcal O\left((\log n_i)(k+1)^{13/4}\right). 
\end{equation}
Thus, this choice yields only polynomial growth for \(\Lambda(P_A)\) as well
\begin{equation}
    \Lambda (P_A) \in \mathcal O \left(|A|^{17/4} \prod_{i=1}^m \log n_i\right).
\end{equation}

How the classic Kronecker product structure of the \emph{Vandermonde} matrices, appearing in the tensorial case, relate to the setup of downward closed index sets is the scope of the next section. 

\section{Kronecker product structure for downward closed index sets}
\label{sec:kronecker_product_structure}
It is well known that, for the lexicographically ordered rectangular index set
\(
    R = \bigtimes_{i=1}^m \{0, \ldots, n_i\},
\)
and for a separable interpolation scheme \((Q_R, P_R)\), the \emph{Vandermonde} matrix and its inverse admit the Kronecker product representations
\begin{equation}
    \label{eq:kronecker_product}
    \bm{\mathcal V} 
    = 
    \bm V^{(1)} \otimes \cdots \otimes \bm V^{(m)},
    \qquad
    \bm{\mathcal V}^{-1} 
    = 
    \big(\bm V^{(1)}\big)^{-1} \otimes \cdots \otimes \big(\bm V^{(m)}\big)^{-1}.
\end{equation}

In this section, we carry out downward closed analogues of these Kronecker product identities.

\subsection{Matrix structure for total, Euclidean, and maximal degree}

First, we illustrate structural examples of \emph{Vandermonde} with respect to the three-variate \emph{Newton} basis in a generic quasi-tensorial grid from Definition~\ref{def:quasi_tensorial_grid}.

\begin{example}
    Figure~\ref{fig:vandermonde_matrices} shows sparsity patterns, not actual matrix values, meaning non-zero entries are represented as black pixels and zero entries as white pixels. \\

    \begin{figure}[ht]
    \centering
    \begin{subfigure}{0.33\textwidth}
        \centering
        \includegraphics{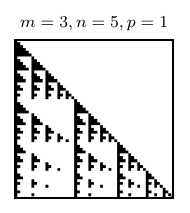}
        \caption{Total degree}
    \end{subfigure}%
    \begin{subfigure}{0.33\textwidth}
        \centering
        \includegraphics{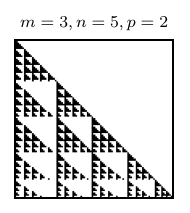}
        \caption{Euclidean degree}
    \end{subfigure}%
    \begin{subfigure}{0.33\textwidth}
        \centering
        \includegraphics{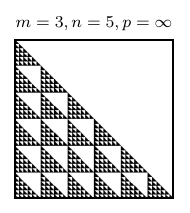}
        \caption{Maximal degree}
    \end{subfigure}%
    \caption{Vandermonde matrices \(\bm {\mathcal V}\) in the three-variate Newton basis of \(p\)-degree \(n \le 5\), shown for three distinct choices of \(p=1,2,\infty\).}
    \label{fig:vandermonde_matrices}
\end{figure}

    The matrices shown in Figure~\ref{fig:vandermonde_matrices} (a) and (b) are submatrices of the tensor-product matrix in Figure~\ref{fig:vandermonde_matrices} (c). Thus, the transition from the tensor-product case \(p=\infty\) to downward closed index sets is directly reflected in the sparsity pattern of the associated \emph{Vandermonde} matrices. For \(p=1\) and \(p=2\), the lower triangular blocks remain visible, but their tails shorten progressively from the north-west to the south-east corner.
\end{example}

This motivates the algebraic construction below, which describes how these submatrices are extracted from the tensor-product \emph{Vandermonde} matrix.

\subsection{Rank embedding and matrix}

The starting point is to describe the embedding of a downward closed index set into a larger one, which can be formalized as follows.

\begin{definition}[Rank embedding]
    \label{def:rank_embedding}
    Let \(A \subset A' \subset \mathbb N_0^m\) be finite, downward closed and both equipped with the same total order \(\preceq\). The \emph{rank embedding} records the position in \(A'\) of each element of \(A\) and is given by
    \begin{equation}
        \label{eq:rank_embedding}
        \varphi_{A, A'} : \{1, \ldots, |A|\} \longrightarrow \{1, \ldots, |A'|\}, 
        \qquad
        j \longmapsto \operatorname{rank}_{A'} (\operatorname{rank}_A^{-1}(j))
    \end{equation}
\end{definition}

To express the rank embedding as a matrix operation, we consider, for \(r \leq s\), the specific binary \emph{Stiefel manifold}
\begin{equation}
    \label{eq:stiefel_manifold}
    \mathrm{St}(r,s) 
    \coloneqq 
    \{\Phi \in \{0,1\}^{r\times s} \mid \Phi \Phi^\top = \bm I^{(r)}\},
\end{equation}
whose elements we call \emph{selection matrices}. The rank embedding \(\varphi_{A, A'}\) then induces a linear map
\begin{equation}
    \varphi^* : \mathbb R^{|A'|} \to \mathbb R^{|A|},
\end{equation}
represented by a specific selection matrix \(\Phi \in \mathrm{St}(|A|,|A'|)\), which is defined below.

\begin{definition}[Rank matrix]
    \label{def:rank_matrix}
    Let the assumptions of Definition~\ref{def:rank_embedding} be fulfilled. 
    We define \(\varphi^*\) as the operation that selects the coordinates of \(\bm v \in \mathbb R^{|A'|}\) according to the image of \(\varphi\)
    \begin{equation}
        \varphi^*(\bm v) 
        \coloneqq 
        \big(v_{\varphi(1)}, \ldots,v_{\varphi(|A|)}\big) \in \mathbb R^{|A|}.
    \end{equation}
    We call the matrix representation of \(\varphi^*\) the \emph{rank matrix} and denote it by
    \begin{equation}
        \Phi_{A, A'} \in \operatorname{St}(|A|, |A'|).
    \end{equation}
    In particular, for the rectangular envelope
    \begin{equation}
        R_A \coloneqq \bigtimes_{i=1}^m \{0, \ldots, n_i\},
        \qquad
        n_i = \max_{\alpha \in A} \alpha_i,
    \end{equation}
    we write 
    \begin{equation} 
        \Phi_A \coloneqq \Phi_{A,R_A}.
    \end{equation}
    \emph{When the index set \(A\) is clear from the context, we simply write \(\Phi\)}.
\end{definition}
As a consequence, given a tensor \(\bm {\mathcal M} \in \mathbb R^{|A'|\times |A'|}\),  the two-sided action of \(\Phi_A\) selects the sub-tensor
\begin{equation}
    \Phi \, \bm {\mathcal M} \, \Phi^\top 
    = 
    \big(
        m_{\,\varphi(j), \varphi(k)}
    \big)_{j, k = 1, \ldots, |A|} \in \mathbb R^{|A|\times|A|}.
\end{equation}

\subsection{Vandermonde matrix}

With the ingredients above we state: 

\begin{proposition}
    \label{prop:selected_kronecker}
    Let \(A \subset \mathbb N_0^m\) be finite, downward closed, and lexicographically ordered. Further, let \((Q_A, P_A)\) be a separable interpolation scheme with \(\bm{\mathcal V}\) denoting its \emph{Vandermonde} matrix, and \(\bm V^{(i)}\) the corresponding univariate Vandermonde matrices. Then, the following identity holds
    \begin{equation}
        \label{eq:selected_tensor_product}
        \bm{\mathcal V}  = \Phi_A \big(\bm V^{(1)} \otimes \cdots \otimes \bm V^{(m)}\big)\Phi_A^\top.
    \end{equation}
\end{proposition}
\begin{proof}
    Let \(\varphi = \varphi_{A,R}\) be the rank embedding from \eqref{eq:rank_embedding}, mapping the rank of each \(\bm \alpha \in A\) with respect to \(A\) to its rank with respect to \(\bm R\). Then, by separability 
    \begin{equation}
        Q_{\bm \beta} (\bm p_{\bm \alpha}) 
        = 
        \prod_{i=1}^m Q_{\beta_i, i} (p_{\alpha_i,i}) 
        = 
        \prod_{i=1}^m v_{\alpha_i+1, \beta_i+1}^{(i)} 
        = 
        (\bm V^{(1)} \otimes \cdots \otimes \bm V^{(m)})_{\varphi(j), \varphi(k)},
    \end{equation}
    where \(\bm \alpha, \bm \beta \in A\) correspond to the ranks \(j, k \in \{1,\ldots, |A|\}\). Thus, Definition~\ref{def:rank_matrix} implies the assertion.
\end{proof}

Equation~\ref{eq:selected_tensor_product} is a selected tensor-product identity. In the \emph{Newton} case, lower triangularity and downward closedness allow the selections \(\Phi_A^\top \Phi_A\) to be inserted between mode-wise tensor factors, leading to the following dimension-wise factorization.

\begin{proposition}
    \label{prop:mode_factorization}
    Let \(A \subset \mathbb N_0^m\) be finite, downward closed, and lexicographically ordered. Further, let \((Q_A, P_A)\) be a \emph{Newton} separable interpolation scheme with \(\bm{\mathcal V}\) denoting its \emph{Vandermonde} matrix, and \(\bm V^{(i)}\) the corresponding univariate \emph{Vandermonde} matrices. We denote the mode-wise tensor factors on the rectangular envelope \(R_A\)
    \begin{equation}
        \bm {\mathcal M}^{(i)}
        \coloneqq
        \bm I^{(n_1+1)}
        \otimes \cdots \otimes
        \bm I^{(n_{i-1}+1)}
        \otimes 
        \bm V^{(i)}
        \otimes
        \bm I^{(n_{i+1}+1)}
        \otimes \cdots \otimes
        \bm I^{(n_m+1)}
    \end{equation}
    Then, the following identity holds
    \begin{equation}
        \bm{\mathcal V}
        =
        \left(
            \Phi_A \, \bm {\mathcal M}^{(1)} \, \Phi_A^\top
        \right)
        \cdots
        \left(
            \Phi_A \, \bm {\mathcal M}^{(m)} \, \Phi_A^\top
        \right).
    \end{equation}
\end{proposition}
\begin{proof}     
    By the tensor-product structure of \(\bm{\mathcal M}^{(i)}\), and the fact that \(\bm V^{(i)}\) is lower triangular, we obtain 
    \begin{equation}
        m^{(i)}_{\bm\alpha,\bm\beta}\neq 0 
        \ \implies \
        \bm\beta\leq \bm\alpha,
    \end{equation}
    Since \(A\) is downward closed, it follows that
    \begin{equation}
        \bm \alpha\in A,\ \bm \beta \notin A \implies m^{(i)}_{\bm\alpha,\bm\beta} = 0,
    \end{equation}
    which yields the identity
    \begin{equation}
        \label{eq:mode_identity}
        \Phi \bm{\mathcal M}^{(i)}
        =
        \Phi \bm{\mathcal M}^{(i)} \Phi^\top \Phi.
    \end{equation}
    Using Proposition~\ref{prop:selected_kronecker} and applying identity~\eqref{eq:mode_identity} successively proves the claim
    \begin{equation}
        \label{eq:selection_procedure}
        \begin{aligned}
            \bm{\mathcal V}
            &=
            \Phi \, \bm{\mathcal M}^{(1)} \cdots \bm{\mathcal M}^{(m)} \, \Phi^\top \\
            &=
            \Phi \, \bm{\mathcal M}^{(1)} \, \Phi^\top \, \Phi \bm{\mathcal M}^{(2)} \cdots \bm{\mathcal M}^{(m)} \, \Phi^\top  \\
            &=
            \left(
                \Phi \bm{\mathcal M}^{(1)}\Phi^\top
            \right)
            \left(
                \Phi \, \bm{\mathcal M}^{(2)} \, \Phi^\top
            \right)
            \cdots
            \left(
                \Phi \, \bm{\mathcal M}^{(m)} \, \Phi^\top
            \right).
        \end{aligned}
    \end{equation}
\end{proof}

\begin{remark}
    For other separable polynomial bases, the procedure in \eqref{eq:selection_procedure} applies \emph{only after} rewriting the corresponding univariate \emph{Vandermonde} matrices through their LU decompositions, see Theorem~\ref{theo:selected_kronecker_lu}. In this case, the computation again passes through the \emph{Newton} basis.
\end{remark}

\subsection{Inverse Vandermonde matrix}

Before carrying out an identity for \(\bm {\mathcal V}^{-1}\), we recall the behavior of inversion under leading principal truncation for triangular matrices.

\begin{proposition} 
    \label{prop:leading_principal_truncation}
    Let \(n, r\in \mathbb N\), \(r \leq n\), \(\Phi = (\bm I^{(r)}, \bm 0) \in \operatorname{St}(r,n)\), and \(\bm V \in \mathbb R^{n\times n}\) be an invertible lower (or upper) triangular matrix. Then, the following identity holds
    \begin{equation}
        \label{eq:leading_principal_truncation}
        \big(\Phi \bm V \Phi^\top\big)^{-1} = \Phi \bm V^{-1} \Phi^\top.
    \end{equation}
\end{proposition}
\begin{proof}
    The matrix \(\Phi \bm V \Phi^\top\) is the leading \(r\times r\) principal submatrix of \(\bm V\). Since \(\bm V\) is invertible and lower triangular, this submatrix is again invertible and lower triangular. By standard block-inversion for triangular matrices, the leading \(r \times r\) principal submatrix of \(\bm V^{-1}\) is precisely \((\Phi \bm V \Phi^\top)^{-1}\). The upper triangular case is analogous. 
\end{proof}

The result below presents the downward closed analogue of the inverse tensor-product identity in \eqref{eq:kronecker_product} for \emph{separable Newton interpolation schemes}. 
It relies on the commutation identity~\eqref{eq:leading_principal_truncation}, which applies to all univariate \emph{Vandermonde} matrices \(\bm V^{(i)}\). 
Without this identity, the result \emph{would fail}.

\begin{proposition}
    \label{prop:selected_kronecker_inverse}
    Let \(A \subset \mathbb N_0^m\) be finite, downward closed, and lexicographically ordered. Further, let \((Q_A, P_A)\) be a separable \emph{Newton} interpolation scheme with \(\bm{\mathcal V}\) denoting its \emph{Vandermonde} matrix, and \(\bm V^{(i)}\) the corresponding univariate \emph{Vandermonde} matrices. Then, the following identity holds
    \begin{equation}
        \label{eq:selected_kronecker_inverse}
        \bm{\mathcal V}^{-1} 
        = 
        \Phi_A \big(
            \big(\bm V^{(1)}\big)^{-1} \otimes \cdots \otimes \big(\bm V^{(m)}\big)^{-1}
        \big) \Phi_A^\top.
    \end{equation}
\end{proposition}
\begin{proof}
    By means of Proposition \ref{prop:selected_kronecker}, we demonstrate \eqref{eq:selected_kronecker_inverse} by showing that
    \begin{equation}
        \big(\bm{\mathcal V}^{-1} \big)_{\bm\alpha,\bm\beta}
        = 
        \prod_{i=1}^m \big(\big(\bm V^{(i)}\big)^{-1}\big)_{\alpha_i+1,\beta_i+1},
        \qquad
        \bm \alpha, \bm \beta \in A.
    \end{equation}
    We consider the separable ansatz
    \begin{equation}
        m_{\bm\alpha,\bm\beta} \coloneqq \prod_{i=1}^m m^{(i)}_{\alpha_i+1,\beta_i+1}, \qquad m^{(i)}_{\alpha_i+1,\beta_i+1} \in \mathbb R,
        \qquad
        \bm \alpha, \bm \beta \in A,
    \end{equation}
    and show that the resulting matrix \(\bm{\mathcal M} = (m_{\bm\alpha,\bm\beta})_{\bm\alpha, \bm\beta \in A} \in \mathbb R^{|A|\times|A|}\) satisfies
    \(\bm{\mathcal V} \bm{\mathcal M} = \bm I\). 
    Since \(Q_A\) is separable, representation \eqref{eq:separable_polynomial_basis} holds for some \(Q_{j,i}\). Thus, we obtain for all \(\bm \alpha, \bm \beta \in A\)
    \begin{equation}
        \sum_{\bm\gamma \in A} Q_{\bm\gamma}(\bm p_{\bm\alpha}) m_{\bm\gamma,\bm\beta} 
        = 
        \sum_{\bm\gamma \in A} \left(\prod_{i=1}^m Q_{\gamma_i,i}(p_{\alpha_i,i})\right) m_{\bm\gamma,\bm\beta}
        =
        \sum_{\bm\gamma \in A} \prod_{i=1}^m v^{(i)}_{ \alpha_i+1, \gamma_i+1} m^{(i)}_{\gamma_i+1,\beta_i+1}.
    \end{equation}
    Grouping this sum with respect to the last coordinate yields
    \begin{equation}
        \sum_{\bm\gamma \in A: \gamma_m=0} \prod_{i=1}^{m-1} v^{(i)}_{\alpha_i+1, \gamma_i+1} m^{(i)}_{\gamma_i+1,\beta_i+1} \sum_{l=1}^{t^{(m)}_{\bm\gamma}} v^{(m)}_{\alpha_m+1, l} m^{(m)}_{l,\beta_m+1},
    \end{equation}
    where
    \begin{equation}
        t^{(i)}_{\bm\gamma} 
        \coloneqq 
        \big|\{j \in \mathbb N_0 \mid (\gamma_1, \ldots, \gamma_{i-1}, j,0,\ldots,0) \in A\}\big|, 
        \qquad 
        \bm\gamma \in A.
    \end{equation}
    Repeating this argument coordinate-wise shows that it is sufficient to have
    \begin{equation}
        \delta_{j,k} 
        = 
        \sum_{l=1}^{t^{(i)}_{\bm\gamma}} v^{(i)}_{j,l} m^{(i)}_{l,k}, \quad j, k = 1, \ldots, n_i+1, \qquad 
        \bm \gamma \in A.
    \end{equation}
    Since \(t^{(i)}_{\bm\gamma} \in \{1, \ldots, n_i+1\}\) for \(\bm\gamma \in A\), this condition is equivalent to 
    \begin{equation}
        \label{app:eq:1}
        \bm I^{(r)} 
        = 
        \left(\Phi^{(r,i)} \, \bm V^{(i)} \, \big(\Phi^{(r,i)}\big)^\top\right) \left(\Phi^{(r,i)} \, \bm M^{(i)} \, \big(\Phi^{(r,i)}\big)^\top\right), 
        \qquad 
        r = 1, \ldots, n_i+1,
    \end{equation}
    where
    \begin{equation}
        \bm M^{(i)}
        \coloneqq
        \big(m^{(i)}_{j,k}\big)_{j,k=1, \ldots, n_i+1},
        \qquad
        \Phi^{(r,i)} \coloneqq (\bm I^{(r)}, \bm0) \in \operatorname{St}(r, n_i+1).
    \end{equation}
    Now, we choose 
    \(
        \bm M^{(i)} = \big(\bm V^{(i)}\big)^{-1}.
    \)
    Since \((Q_A, P_A)\) is a separable \emph{Newton} interpolation scheme, all \(\bm V^{(i)}\) are lower triangular matrices. 
    Hence, by Proposition~\ref{prop:leading_principal_truncation}
    \begin{equation}
        \Phi^{(r,i)} \, \big(\bm V^{(i)}\big)^{-1} \, \big(\Phi^{(r,i)}\big)^\top 
        = 
        \big(\Phi^{(r,i)} \, \bm V^{(i)}\,  \big(\Phi^{(r,i)}\big)^\top \big)^{-1}.
    \end{equation}
    Thus \eqref{app:eq:1} holds, and consequently 
    \(
        \bm{\mathcal V} \, \bm{\mathcal M} = \bm I,
    \)
    proving \eqref{eq:selected_kronecker_inverse}.
\end{proof}

\begin{remark}
    The identity for the inverse \emph{Vandermonde} matrix \(\bm{\mathcal V}^{-1}\) in \eqref{eq:selected_kronecker_inverse} admits, analogously to Proposition~\ref{prop:mode_factorization}, the mode-wise factorization
    \begin{equation}
        \bm{\mathcal V}^{-1} 
        = 
        \left(
            \Phi_A 
            \, 
            \big(\bm {\mathcal M}^{(1)}\big)^{-1} 
            \, 
            \Phi_A^\top
        \right)
        \cdots
        \left(
            \Phi_A 
            \, 
            \big(\bm {\mathcal M}^{(m)}\big)^{-1} 
            \, 
            \Phi_A^\top\right).
    \end{equation}
\end{remark}

\subsection{LU-based formulas for separable schemes}

The \emph{fast Newton transform} relies on the hierarchical triangular structure of the \emph{Vandermonde} matrix, which is not directly available for general separable interpolation schemes.
However, univariate LU decompositions can be used to reduce the general case to the separable \emph{Newton} setting required by the \emph{FNT}.
This yields the following two Kronecker product formulas, at the cost of performing \(m\) univariate LU transforms as described already in Section~\ref{subsec:main_results}.

\begin{theorem}
    \label{theo:selected_kronecker_lu}
    Let \(A \subset \mathbb N_0^m\) be finite, downward closed, and lexicographically ordered. 
    Further, let \((Q_A, P_A)\) be a separable interpolation scheme with \(\bm{\mathcal V}\) denoting its \emph{Vandermonde} matrix, and  \(\bm V^{(i)}\) the corresponding univariate \emph{Vandermonde} matrices. Together with the LU decompositions
    \begin{equation}
        \bm V^{(i)} = \bm L^{(i)} \bm U^{(i)}, 
        \qquad 
        \bm L^{(i)}, \bm U^{(i)} \in \mathbb R^{(n_i+1)\times(n_i+1)}.
    \end{equation}
    the following identities hold:
    \begin{enumerate}[label=(\roman*)]
        \item \(
            \bm{\mathcal V} 
            = 
            \Phi_A \left(
                \bm L^{(1)} 
                \otimes \cdots \otimes 
                \bm L^{(m)}
            \right) \Phi_A^\top 
            \,
            \Phi_A \left(
                \bm U^{(1)} 
                \otimes \cdots \otimes
                \bm U^{(m)}
            \right) \Phi_A^\top\),
        \item \(
            \bm{\mathcal V}^{-1} 
            = 
            \Phi_A \left(
                \big(\bm U^{(1)}\big)^{-1} 
                \otimes \cdots \otimes 
                \big(\bm U^{(m)}\big)^{-1}
            \right) \Phi_A^\top 
            \,
            \Phi_A \left(
                \big(\bm L^{(1)}\big)^{-1} 
                \otimes \cdots \otimes 
                \big(\bm L^{(m)}\big)^{-1}
            \right) \Phi_A^\top\).
    \end{enumerate}
\end{theorem}
\begin{proof}
    Defining
    \begin{equation}
        \bm{\mathcal L} 
        \coloneqq 
        \bm L^{(1)} \otimes \cdots \otimes \bm L^{(m)},
        \qquad
        \bm{\mathcal U} 
        \coloneqq 
        \bm U^{(1)} \otimes \cdots \otimes \bm U^{(m)},
    \end{equation}
    and applying Proposition~\ref{prop:selected_kronecker} yields
    \begin{equation}
        \bm{\mathcal V} 
        = 
        \Phi \left(
            \bm V^{(1)} 
            \otimes \cdots \otimes 
            \bm V^{(m)}
        \right) \Phi^\top 
        = 
        \Phi \left(
            \big(\bm L^{(1)} \bm U^{(1)}\big)
            \otimes \cdots \otimes
            \big(\bm L^{(m)} \bm U^{(m)}\big)
        \right) \Phi^\top 
        = 
        \Phi \, \bm{\mathcal L} \, \bm{\mathcal U} \, \Phi^\top.
    \end{equation}
    To prove \emph{(i)}, it remains to justify that the selection \(\Phi^\top \Phi\) may be inserted between \(\bm {\mathcal L}\) and \(\bm {\mathcal U}\). Let \(\bm\alpha, \bm\beta \in A\), then
    \begin{equation}
        (\Phi \, \bm{\mathcal L} \, \bm{\mathcal U} \, \Phi^\top)_{\bm\alpha, \bm\beta} 
        = 
        \sum_{\bm\gamma \in R_A} l_{\bm \alpha, \bm \gamma} u_{\bm \gamma, \bm \beta}.
    \end{equation}
    Since \(\bm{\mathcal L}\) is a Kronecker product of lower triangular matrices and \(A\) is downward closed, we obtain
    \begin{equation}
        l_{\bm \alpha, \bm \gamma} \neq 0 \implies \bm \gamma \le \bm \alpha \implies \bm \gamma \in A.
    \end{equation}
    Hence, the sum over \(R_A\) can be restricted to \(A\)
    \begin{equation}
        (\Phi \, \bm{\mathcal L} \, \bm{\mathcal U} \, \Phi^\top)_{\bm\alpha, \bm\beta} 
        = 
        \sum_{\bm\gamma \in A} l_{\bm \alpha, \bm \gamma} u_{\bm \gamma, \bm \beta} 
        = 
        (\Phi \, \bm{\mathcal L} \, \Phi^\top \, \Phi \, \bm{\mathcal U} \, \Phi^\top)_{\bm\alpha, \bm\beta},
    \end{equation}
    and, equivalently
    \begin{equation}
        \bm{\mathcal{V}} 
        = 
        \left(\Phi \, \bm{\mathcal L} \, \Phi^\top\right) \left(\Phi \, \bm{\mathcal U} \, \Phi^\top\right).
    \end{equation}
    Finally, combining \emph{(i)} with Proposition~\ref{prop:selected_kronecker_inverse} yields \emph{(ii)}
    \begin{equation}
        \bm{\mathcal V}^{-1} 
        = 
        \left(
            \Phi \, \bm{\mathcal L} \, \Phi^\top \, \Phi \, \bm{\mathcal U} \, \Phi^\top
        \right)^{-1} 
        = 
        \left(
            \Phi \, \bm{\mathcal U} \, \Phi^\top\right)^{-1} \left(\Phi \, \bm{\mathcal L} \, \Phi^\top\right)^{-1} 
        = 
        \Phi \, \bm{\mathcal U}^{-1} \, \Phi^\top \, \Phi \bm{\mathcal L}^{-1} \, \Phi^\top.
    \end{equation}
\end{proof}

The generalization of the Kronecker product given herein suggests the following construction of the \emph{fast Newton transform}.

\section{Construction of the fast Newton transform}
\label{sec:construction_of_the_fnt}

In this section, we prove the main complexity results stated in Theorem~\ref{theo:fnt} and Theorem~\ref{theo:fastdiff}. 

We begin with the proof of Theorem~\ref{theo:fnt}. The idea is to rewrite the \emph{Newton} transform and its inverse as sequences of mode-wise matrix applications.
We then show that suitable index permutations reduce each of these operations to a last-mode matrix application.

\fastnewtontransform*
\begin{proof}
    We reformulate the computational challenges \ref{C1} and \ref{C2}, as the computation of \(\bm c\) and \(\bm f\), respectively, in the form
    \begin{equation}
        \label{eq:reformulation_of_computational_challenges}
        \bm c 
        =
        \Phi_A\left(
            \big(\bm V^{(1)}\big)^{-1} \otimes \cdots \otimes \big(\bm V^{(m)}\big)^{-1}
        \right)\Phi_A^\top \bm f,
        \qquad 
        \bm f 
        = 
        \Phi_A\left(
            \bm V^{(1)} \otimes \cdots \otimes \bm V^{(m)}
        \right)\Phi_A^\top \bm c.
    \end{equation}
    For \emph{Newton} interpolation, Proposition~\ref{prop:mode_factorization} allows these selected tensor products to be decomposed into products of mode-wise factors.
    Starting with \(\bm f_m = \bm f\) and \(\bm c_m = \bm c\), we obtain the iterative schemes
    \begin{equation}
        \label{eq:iterative_scheme}
        \begin{aligned}
            \bm f_{i-1} 
            &\coloneqq
            \Phi \left(
                \bm I^{(n_1+1)}
                \otimes \cdots \otimes
                \bm I^{(n_{i-1}+1)}
                \otimes 
                \big(\bm V^{(i)}\big)^{-1}
                \otimes
                \bm I^{(n_{i+1}+1)}
                \otimes \cdots \otimes
                \bm I^{(n_m+1)}
            \right)\Phi^\top \bm f_i, \\
            \bm c_{i-1}
            &\coloneqq
            \Phi\left(
                \bm I^{(n_1+1)}
                \otimes \cdots \otimes
                \bm I^{(n_{i-1}+1)}
                \otimes 
                \bm V^{(i)}
                \otimes
                \bm I^{(n_{i+1}+1)}
                \otimes \cdots \otimes
                \bm I^{(n_m+1)}
            \right)\Phi^\top \bm c_i, \\
        \end{aligned}
    \end{equation}
    delivering the \emph{Newton} coefficients \(\bm c = \bm f_0\) and the function evaluations \(\bm f = \bm c_0\). 
    Further, we use index permutations \(\psi_i\) to move the \(i\)th coordinate direction to the last position. 
    Let \(\Psi^{(i)}\) denote the corresponding permutation matrix on \(A\).
    Then, the previous schemes can be written equivalently as
     \begin{equation}
        \begin{aligned}
            \bm f_{i-1} 
            &= 
            \big(\Phi \Psi^{(i)}\big)
            \big(
                \bm I^{(n_{i+1}+1)}
                \otimes \cdots \otimes
                \bm I^{(n_m+1)}
                \otimes 
                \bm I^{(n_1+1)}
                \otimes \cdots \otimes
                \bm I^{(n_{i-1}+1)}
                \otimes
                \big(\bm V^{(i)}\big)^{-1}
            \big)
            \big(\Phi\Psi^{(i)}\big)^\top \bm f_i, \\
            \bm c_{i-1}
            &= 
            \big(\Phi \Psi^{(i)}\big)
            \big(
                \bm I^{(n_{i+1}+1)}
                \otimes \cdots \otimes
                \bm I^{(n_m+1)}
                \otimes 
                \bm I^{(n_1+1)}
                \otimes \cdots \otimes
                \bm I^{(n_{i-1}+1)}
                \otimes
                \bm V^{(i)}
            \big)
            \big(\Phi \Psi^{(i)}\big)^\top \bm c_i.
        \end{aligned}
    \end{equation}
    By Proposition~\ref{prop:index_permutations}, the index permutations \(\psi_1, \ldots, \psi_m\) can be precomputed in \(\mathcal O(m|A|)\) operations and require \(\mathcal O (m|A|)\) storage. Applying each such permutation costs \(\mathcal O (|A|)\) steps. Moreover, by Proposition~\ref{prop:last_mode_matrix}, the last-mode application of \(\bm V^{(i)}\), respectively, \(\big(\bm V^{(i)}\big)^{-1}\), costs \(\mathcal O (n_i |A|)\). Summing over all modes \(i=1,\ldots,m\) gives
    \begin{equation}
        \mathcal O ((n_1+\cdots+n_m)|A|) = \mathcal O (m \overline n |A|),
    \end{equation}
    which is the claimed time complexity.
\end{proof}

As a result of the proof of Theorem~\ref{theo:fnt}, we obtain the following complexity result for polynomial differentiation.

\fastdiff*
\begin{proof}
    Since the \emph{Newton} basis \(Q_A\) is separable, differentiation with respect to one coordinate acts only on the corresponding univariate factor.
    For each coordinate, let 
    \(
        \bm D^{(i)} = \bigl(d_{l,k}^{(i)}\bigr)_{l,k=0,\ldots,n_i} \in \mathbb R^{(n_i+1)\times(n_i+1)}
    \)
    denote the univariate coefficient differentiation matrix in the \emph{Newton} basis, which is defined through
    \begin{equation}
        \frac{\mathrm d}{\mathrm d x} Q_{k,i}(x) = \sum_{l=0}^{n_i} d_{l,k}^{(i)} Q_{l,i}(x).
    \end{equation}
    Since \(\deg Q_{k,i}=k\), \(d_{l,k}^{(i)}=0\) for \(l\ge k\). Hence \(\bm D^{(i)}\) is strictly upper triangular.
    Further, leveraging the degree-graded structure of the \emph{Newton} basis, each matrix \(\bm D^{(i)}\) can be computed in \(\mathcal O(n_i^2)\) operations, as detailed in \cite{Amiraslani2016}.
    The coefficients \(\bm c' \in \mathbb R^{|A|}\) of \(\partial_{x_i} q\) can then be computed by
    \begin{equation}
        \bm c'
        =
        \Phi_A(
            \bm I^{(n_1+1)}
            \otimes \cdots \otimes
            \bm I^{(n_{i-1}+1)}
            \otimes 
            \bm D^{(i)}
            \otimes
            \bm I^{(n_{i+1}+1)}
            \otimes \cdots \otimes
            \bm I^{(n_m+1)}
        )\Phi_A^\top \bm c.
    \end{equation}
    Combining Proposition~\ref{prop:index_permutations} with Proposition~\ref{prop:last_mode_matrix}, analogously to the proof of Theorem~\ref{theo:fnt}, this mode-wise application of \(\bm D^{(i)}\) can be carried out in \(\mathcal O(n_i |A|)\) operations.
\end{proof}

The remainder of this section is devoted to the construction of the index permutations \(\psi_1,\ldots,\psi_m\) and to the last mode matrix application routine, both used in the proofs of Theorem~\ref{theo:fnt} and Theorem~\ref{theo:fastdiff}.

\subsection{Computing the index permutations}
The index permutations \(\psi_1, \ldots, \psi_m\) are induced by cyclic coordinate shifts of downward closed index sets. 
The cyclic coordinate shift that moves the \(i\)th coordinate to the last position is given by
\begin{equation}
    \label{eq:cyclic_coordinate_shift}
    \eta_i: \mathbb N_0^m \longrightarrow 
    \mathbb N_0^m,
    \qquad
    \bm \alpha
    \longmapsto
    (\alpha_{i+1}, \ldots, \alpha_m, \alpha_1, \ldots, \alpha_i).
\end{equation}
Let \(A \subset \mathbb N_0^m\) be finite, downward closed and lexicographically ordered. We denote by 
\begin{equation}
    A^{(i)} \coloneqq \eta_i(A), 
\end{equation}
the corresponding shifted index set \emph{equipped with the lexicographic order}.
The induced index permutation \(\psi_i \in S_{|A|}\) is then determined by satisfying
\begin{equation}
    \label{eq:index_permutation}
    \operatorname{rank}_{A^{(i)}}\!\big(\eta_i(\bm \alpha)\big)
    =
    \psi_i\big(\operatorname{rank}_A(\bm \alpha)\big),
\end{equation}
for all \(\bm \alpha \in A\). Let \(\Psi^{(i)} \in \operatorname{St}(|A|,|A|)\) denote the corresponding permutation matrix, being the transformation matrix of the linear map
\begin{equation}
    \label{eq:permutation_map}
    \psi_i^*(\bm x)
    \coloneqq
    \bigl(x_{\psi_i(j)}\bigr)_{j=1, \ldots, |A|},
    \qquad
    \bm x \in \mathbb R^{|A|}.
\end{equation}
\(\Psi^{(i)}\) realizes the rotation of coordinates in \(\mathbb R^m\) that moves the \(i\)th coordinate to the last position. In particular, applying \(\Psi^{(i)}\) or its inverse \((\Psi^{(i)})^{-1}\) to a vector in \(\mathbb R^{|A|}\) requires only \(\mathcal O(|A|)\) operations.

Thus, the remaining task is to construct the permutations \(\psi_i\) efficiently.
This construction is not immediate, since after rotating the coordinates, the set \(\eta_i(A)\) has to be ordered lexicographically again. 
A naive approach would rotate all multi-indices, sort the resulting set, and compare the two orderings to recover \(\psi_i\). Since the sorted objects are multi-indices of length \(m\), this leads to a cost of order
\(
    \mathcal O(m |A| \log |A|)
\)
for each permutation, and hence to
\begin{equation}
    \mathcal O(m^2 |A| \log |A|)
\end{equation}
for constructing all index permutations \(\psi_1,\ldots,\psi_m\). 

Thus, we provide an alternative low cost construction as follows.

\begin{proposition}
    \label{prop:index_permutations}
    Let \(A \subset \mathbb N_0^m\) be finite, downward closed and lexicographically ordered. The permutations
    \(
        \psi_1, \ldots, \psi_m
    \)
    defined through \eqref{eq:index_permutation} can be constructed in \(\mathcal O(m|A|)\) operations.
\end{proposition}
\begin{proof}
    Rather than constructing each permutation \(\psi_i\) directly, we compute the cyclic coordinate shift \(\eta_i\) as a composition of elementary cyclic coordinate shifts. Hence, the induced permutation \(\psi_i\) can be written as
    \begin{equation}
        \label{eq:index_tranpositions}
        \psi_i
        =
        \tau_1 \circ \tau_2 \circ \cdots \circ \tau_{m-i},
        \qquad
        \psi_m \equiv \operatorname{id}.
    \end{equation}
    Here, each \(\tau_j\) denotes the index permutation obtained by moving the coordinate currently in the last position to the first position, equivalently by applying one cyclic right shift, and then restoring the lexicographic ordering of the index set.
    Algorithm~\ref{alg:index_transposition} constructs each such permutation in \(\mathcal O(|A|)\) operations. Thus, all permutations \(\psi_1, \ldots, \psi_m\) are obtained by successive compositions of the  \(\tau_j\) in \(\mathcal O(m|A|)\) operations.
\end{proof}

Next, we compute the index permutations \(\tau_1, \ldots, \tau_m\).

\subsection{Computing the elementary index permutations}

In analogy with tensor analysis~\cite{Kolda2009}, where tubes describe
variation along the last mode of a tensor, we introduce \emph{tube size projections} for downward closed index sets \(A \subset \mathbb N_0^m\).
The tubes of \(A\) are obtained by fixing all coordinates except the last one.
Equivalently, they are the set-theoretic fibers of the \emph{coordinate projection}
\begin{equation}
    \label{eq:coordinate_projection}
    \pi_m: \mathbb N_0^m \rightarrow \mathbb N_0^{m-1}, 
    \qquad 
    \bm \alpha \longmapsto (\alpha_1, \ldots, \alpha_{m-1}).
\end{equation}
We formalize this partition and record the resulting tube sizes as ordered integer compositions~\cite{Stanley2011}.
\begin{definition}[Tube size projections]
    \label{def:tube_size_projections}
    Let \(A \subset \mathbb N_0^m\) be finite, downward closed, and lexicographically ordered. We define the \emph{tube size projection} of \(A\) by
    \begin{equation}
        \bm t(A) \coloneqq \big(t_{\bm\beta}\big)_{\bm \beta \in \pi_m (A)}, 
        \qquad 
        t_{\bm\beta} \coloneqq |\{\bm \alpha \in A \ : \ \pi_m (\bm \alpha) = \bm \beta \}|,
    \end{equation}
    with \(\pi_m(A)\) being lexicographically ordered. We omit the dependence of \(A\) when it is clear from the context.
\end{definition}
Before turning to the construction of the tube size projections, we first show how they are used. Assuming that these projections are given, Algorithm~\ref{alg:index_transposition} computes the elementary index permutations \(\tau_j\) by following a standard bucket-sort principle~\cite{Knuth1997}.

\vspace{0.5em}
    
\begin{algorithm}[ht]
    \small
    \caption{Computing the elementary index permutation $\tau_j$}
    \label{alg:index_transposition}
    \KwIn{The cyclically shifted index set $A^{(j)} = \eta_j(A)$}
    \KwOut{The elementary index permutation $\tau_j$}
    $N  \gets |A^{(j)}|$\;
    $\bm t(A^{(j)}) \gets (t_1,\ldots, t_N)$\;
    $M \gets \max_{l=1,\ldots,N} \, t_l$\;
    Precompute offsets  $\bm s \gets (s_1,\ldots,s_N)$ with $s_l \gets \sum_{r=1}^{l-1} t_r$\;
    Construct buckets $B_k \gets \{l \in \{1,\ldots,N\} : f_l = k\}$ for $k = 1, \ldots, M$\;
    Initialize an ordered active list $L_1 \gets (1, 2, \ldots, N)$\;
    $r \gets 1$\;
    \For{$k = 1,\ldots,M$}{
        \For{$l \in L_k$}{
            $\tau_j(r) \gets s_l + k$\;
            $r \gets r+1$\;
        }
        $L_{k+1} \gets L_k \setminus B_k$\;
    }
    \Return $\tau_j$\;
\end{algorithm}

\vspace{0.5em}

In Algorithm~\ref{alg:index_transposition}, \(A^{(j)}\) is assumed to be given in lexicographic order.
In practice, this order is obtained recursively from \(A^{(j-1)}\) by applying the corresponding cyclic coordinate shift and the previously constructed elementary index permutation \(\tau_{j-1}\). 
Since both steps require only a single traversal of the index set, the transition from \(A^{(j-1)}\) to \(A^{(j)}\) is performed in \(\mathcal O(|A|)\).

This yields the following time complexity result.

\begin{proposition}
    Let \(A \subset \mathbb N_0^m\) be finite, downward closed, and lexicographically ordered. Given the tube size projection \(\bm t(A^{(j)})\), Algorithm~\ref{alg:index_transposition} computes the elementary permutation \(\tau_j\) from~\eqref{eq:index_tranpositions} in \(\mathcal O(|A|)\) operations.
\end{proposition}
\begin{proof}
    Algorithm~\ref{alg:index_transposition} constructs \(\tau_j\) from the tube size projection
    \begin{equation}
        \bm t(A^{(j)}) = (t_1, \ldots, t_N).
    \end{equation}
    The offsets \(s_l\) mark the first position of the \(l\)th fiber in the lexicographic ordering of \(A^{(j)}\). Then, for every pair \((k,l)\) satisfying
    \(
        1 \le k \le t_l,
    \)
    the index \(s_l + k\) is the position of the \(k\)th element in the \(l\)th fiber.
    The active list \(L_k\) is maintained such that, during the \(k\)th outer iteration
    \begin{equation}
        L_k = \{\, l \in \{1,\ldots,N\} : t_l \ge k \,\},
    \end{equation}
    Thus, the inner loop writes precisely all indices \(s_l + k\) with \(t_l \ge k\). After this iteration, all indices \(l\) with \(t_l = k\) are removed, since they do not satisfy \(t_l \ge k + 1\).
    Consequently, the number of operations for \(\tau_j\) is
    \begin{equation}
        \sum_{k=1}^M |L_k|
        =
        \sum_{k=1}^M
        \big|\{l \in \{1,\ldots,N\} : t_l \ge k\}\big|
        =
        \sum_{l=1}^N t_l
        =
        |A|.
    \end{equation}
    Moreover, each fiber index is inserted into exactly one bucket and removed from the active lists at most once. Assuming constant time removal from the active lists, the time complexity is
    \begin{equation}
        \mathcal O\left(
            |A| + N + M
        \right)
        =
        \mathcal O(|A|),
    \end{equation}
    since \(N \le |A|\), \(M \le |A|\), and \(\sum_{k=1}^M |L_k| = |A|\).
\end{proof}

It remains to show how the tube size projections can be computed. 

\subsection{Computing the tube size projection}

For finite downward closed index sets \(A \subset \mathbb N_0^m\), the tube size projection \(\bm t(A)\) can be obtained by Algorithm~\ref{alg:tube_projection}.

\vspace{0.5em}
    
\begin{algorithm}[ht]
    \small
    \caption{Computing the tube size projection $\bm t(A)$}
    \label{alg:tube_projection}
    \KwIn{
        finite non-empty downward closed and lexicographically ordered index set \(A \subset \mathbbm N_0^m\)
    }
    \KwResult{
        Tube size projection \(\bm t(A)\)
    }
    $\bm t \gets \bm 0 \in \mathbbm N^{|A|}$\;
    $r \gets 1$\;
    $j \gets 1$\;
    \For{$\bm\alpha \in A \setminus \{\bm 0\}$}{
        \If{$\alpha_m > 0$}{
            $r \gets r+1$\;
        }
        \Else{
            $t_j \gets r$\;
            $j \gets j+1$\;
            $r \gets 1$\;
        }
    }
    $t_j \gets r$\;
    \Return $(t_1,\ldots,t_j)$\;
\end{algorithm}

\vspace{0.5em}

\begin{proposition}
    \label{prop:tube_projection_complexity}
    Let \(A \subset \mathbb N_0^m\) be finite, downward closed, and
    lexicographically ordered. Then, Algorithm~\ref{alg:tube_projection}
    computes the tube-size projection \(\bm t(A)\) in \(\mathcal O(|A|)\)
    steps.
\end{proposition}
\begin{proof}
    Since \(A\) is lexicographically ordered, all indices with the same first
    \(m-1\) coordinates appear consecutively. Each tube with respect to
    the last coordinate forms one contiguous block. Since \(A\) is downward
    closed, each tube contains its minimal element, whose last
    coordinate is zero. Thus, a new tube starts exactly when an index \(\bm \alpha \in A\) with \(\alpha_m = 0\) is encountered.

    Algorithm~\ref{alg:tube_projection} traverses \(A\) once, starts a new count at
    each such index, and appends the length of the previous block to
    \(\bm t(A)\). Therefore, the total cost is \(\mathcal O(|A|)\).
\end{proof}

\begin{remark}
    For the particular choice \(A=A_{m,n,p}\), the construction simplifies further. 
    Since \(A_{m,n,p}\) is invariant under coordinate permutations, all cyclically shifted index sets have the same tube size projection. Consequently, the elementary index permutations \(\tau_j\) constructed from these projections coincide. Hence, the permutations \(\psi_i\) can be obtained by repeated composition of one fixed elementary permutation.
\end{remark}

With the construction of \(\tau_1,\ldots,\tau_m\) in place, we now assemble the routine for applying matrices along the last mode.

\subsection{Last mode matrix application}

We present the main routine in Algorithm~\ref{alg:last_mode_matrix}, which applies a matrix along the coordinate tubes of a lexicographically ordered downward closed index set.

\vspace{0.5em}
    
\begin{algorithm}[ht]
    \small
    \caption{Last mode matrix application}
    \label{alg:last_mode_matrix}
    \KwIn{
        Matrix $\bm V \in \mathbb R^{n\times n}$; \
        Vector $\bm x \in \mathbbm R^{|A|}$; \
        Tube size projection $\bm t(A) = (t_1,\ldots,t_q)$; \
    }
    \KwResult{
        Last-mode application of $\bm V$ to $\bm x$
    }
    Offsets $s_1,\ldots,s_{q+1}$ with $s_1 \gets 0$ and $s_{l+1} \gets s_l + t_l$;
    $\bm y \gets \bm 0 \in \mathbbm R^{s_{q+1}}$\;
    \For{$l = 1,\ldots,q$}{
        $\bm x_l \gets (x_{s_l},\ldots,x_{s_{l+1}-1})^\top$\;
        Let $\bm V_{t_l} \in \mathbbm R^{t_l \times t_l}$ denote the leading principal submatrix of $\bm V$\;
        $\bm y_l \gets \bm V_{t_l} \, \bm x_l$\;
        $(y_{s_l},\ldots,y_{s_{l+1}-1})^\top \gets \bm y_l$\;
    }
    \Return $\bm y$\;
\end{algorithm}

\vspace{0.5em}

We conclude with the corresponding time complexity result.

\begin{proposition}
    \label{prop:last_mode_matrix}
    Let \(A \subset \mathbb N_0^m\) be finite, downward closed, and lexicographically ordered. Furthermore, let \(\bm V \in \mathbb R^{(n_m+1)\times(n_m+1)}\). Then, for every \(\bm x \in \mathbb R^{|A|}\), the matrix-vector product
    \begin{equation}
        \label{eq:last_mode_product}
        \Phi \big(
            \bm I^{(n_1+1)}
            \otimes \cdots \otimes
            \bm I^{(n_{m-1}+1)}
            \otimes
            \bm V
        \big) \Phi^\top
        \bm x,
    \end{equation}
    can be computed in \(\mathcal O(n_m |A|)\) operations.
\end{proposition}
\begin{proof}
    Since \(A\) is lexicographically ordered, the last-coordinate tubes appear as contiguous blocks. Let
    \begin{equation}
        \bm t(A)=(t_1,\ldots,t_q),
    \end{equation}
    denote the corresponding tube size projection, where \(q=|\pi_m (A)|\). Thus, it holds
    \begin{equation}
        \sum_{l=1}^q t_l = |A|,
        \qquad
        t_l \le n_m + 1.
    \end{equation}
    Accordingly, every vector \(\bm x \in \mathbb R^{|A|}\) decomposes into contiguous subvectors
    \begin{equation}
        \bm x = (\bm x^{(1)}, \ldots, \bm x^{(q)}),
        \qquad
        \bm x^{(l)} \in \mathbb R^{t_l}.
    \end{equation}
    Consider now the embedded vector \(\Phi^\top \bm x \in \mathbb R^{|R|}\). For the \(l\)th tube, the embedding \(\Phi^\top\) extends the block \(\bm x^{(l)}\) to the full last-coordinate fiber by appending zeros. Hence, applying the last-mode factor \(\bm V\) and restricting back with \(\Phi\) gives
    \begin{equation}
        \left[
            \Phi
            \left(
                \bm I^{(n_1+1)}
                \otimes \cdots \otimes
                \bm I^{(n_{m-1}+1)}
                \otimes
                \bm V
            \right)
            \Phi^\top \bm x
        \right]^{(l)}
        =
        \begin{pmatrix}
            \bm I^{(t_l)} & \bm 0
        \end{pmatrix}
        \bm V
        \begin{pmatrix}
            \bm x^{(l)}\\
            \bm 0
        \end{pmatrix}
        =
        \bm V_{t_l}\bm x^{(l)},
    \end{equation}
    where \(\bm V_{t_l}\in\mathbb R^{t_l\times t_l}\) denotes the leading principal submatrix of \(\bm V\).
    Therefore, Algorithm~\ref{alg:last_mode_matrix} computes the product in \eqref{eq:last_mode_product} by applying the matrix-vector products
    \(
        \bm V_{t_l}\bm x^{(l)}
    \)
    on all tubes. Hence, the total number of operations is bounded by
    \begin{equation}
        \mathcal O\left(
            \sum_{l=1}^q t_l^2
        \right).
    \end{equation}
    Using \(t_l \leq n_m+1\), we obtain
    \begin{equation}
        \sum_{l=1}^q t_l^2
        \leq
        (n_m+1)\sum_{l=1}^q t_l
        =
        (n_m+1)|A|.
    \end{equation}
    Therefore, the product in \eqref{eq:last_mode_product} can be computed in \(\mathcal O(n_m|A|)\) operations.
\end{proof}

Given the theoretically proven time complexities, we verify them through numerical experiments.

\section{Numerical experiments}
\label{sec:numerical_experiments}

The Python package \emph{lpFun}~\cite{Hofmann2026lpFun} implements the specifications described in Section~\ref{sec:implementation_choices}.
Throughout the experiments, we fix \(Q_A\) to be the \emph{Newton} basis and \(P_A\) to be generated by the Leja-ordered \emph{Chebyshev}-\emph{Lobatto} nodes.

We assess performance by reporting the maximum relative error against runtime in milliseconds. For each test function, we record two different runtimes. First, we measure the time required to compute the coefficients of the interpolant, corresponding to \ref{C1}. Second, we measure the total time required to compute the coefficients of the derivative of the interpolant with respect to the first coordinate and subsequently evaluate this derivative in the grid \(P_A\), corresponding to \ref{C2}. Both runtimes are averaged over 10 repetitions. We evaluate the maximum relative error over uniformly randomly sampled points and the number of evaluation points is chosen according to the spatial dimension
\begin{equation}
    10^6 \quad \text{for } m=1,2,3,
    \qquad
    10^5 \quad \text{for } m=4,
    \qquad
    10^4 \quad \text{for } m=5,
    \qquad
    10^3 \quad \text{for } m=6.
\end{equation}
For each configuration, the degree is increased until either the maximum relative error of the interpolant no longer improves for two consecutive degrees, or until it reaches \(10^{-15}\).

All experiments were conducted on the same machine equipped with an Apple M2 Pro processor featuring 12 cores (8 performance, 4 efficiency) and 16 GB of unified memory.

\subsection{Comparison of \(\ell^p\)-type index sets}

We first investigate how the choice of the downward closed index set \(A\) affects the numerical performance of the \emph{FNT}. Therefore, we compare the following \(\ell^p\)-type index sets in spatial dimensions \(m = 1, \ldots, 6\)
\begin{equation}
    A_{m,n,p}
    =
    \{
        \bm\alpha \in \mathbb N_0^m
        \ : \
        \|\bm\alpha\|_p \le n
    \},
    \qquad 
    p \in \{1,2,\infty\}.
\end{equation}

We use the following entire test functions:

\begin{description}
    \item[\textbf{A}.]
        Separable sine product function
        \begin{equation}
        f(\bm x)
        =
        \prod_{i=1}^m \sin(\pi x_i).
    \end{equation}
    
    \item[\textbf{B}.]
    Mixed exponential cosine function
    \begin{equation}
        f(\bm x)
        =
        \exp\left(
            \frac{1}{m}\sum_{i=1}^m x_i
        \right) \cos\left(
            \frac{1}{m}\sum_{i=1}^m x_i
        \right).
    \end{equation}
    
    \item[\textbf{C}.]
    Anisotropic exponential function 
    \begin{equation}
        f(\bm x)
        =
        \exp\left(\sum_{i=1}^m \frac{x_i}{i^2}\right).
    \end{equation}
\end{description}

\begin{figure}[!htbp]
    \centering

    \includegraphics[width=0.9\textwidth]{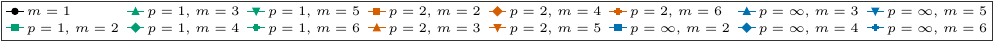}

    \vspace{0.5em}

    \begin{subfigure}{\textwidth}
        \centering
        \includegraphics[width=0.9\textwidth]{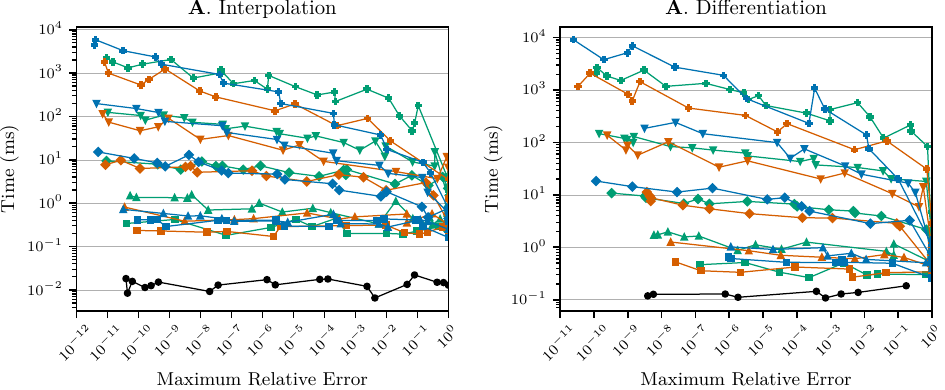}
        \caption{Separable sine product function.}
        \label{fig:results_separable_sine_product}
    \end{subfigure}

    \vspace{0.8em}

    \begin{subfigure}{\textwidth}
        \centering
        \includegraphics[width=0.9\textwidth]{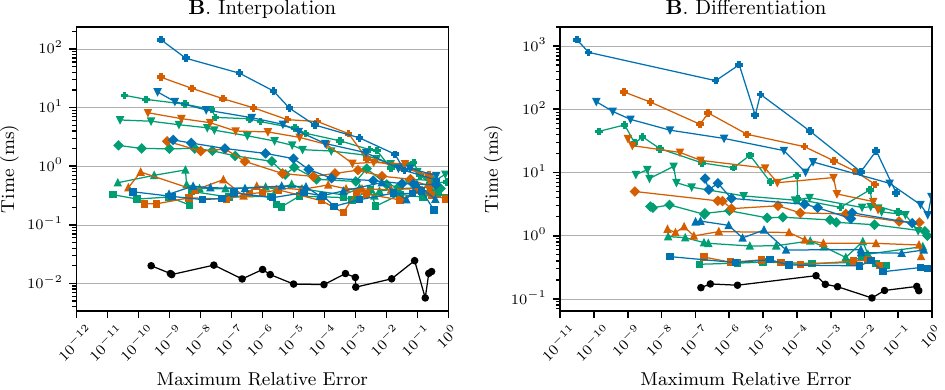}
        \caption{Mixed exponential cosine function.}
        \label{fig:results_hyperbolic_tangent}
    \end{subfigure}

    \vspace{0.8em}

    \begin{subfigure}{\textwidth}
        \centering
        \includegraphics[width=0.9\textwidth]{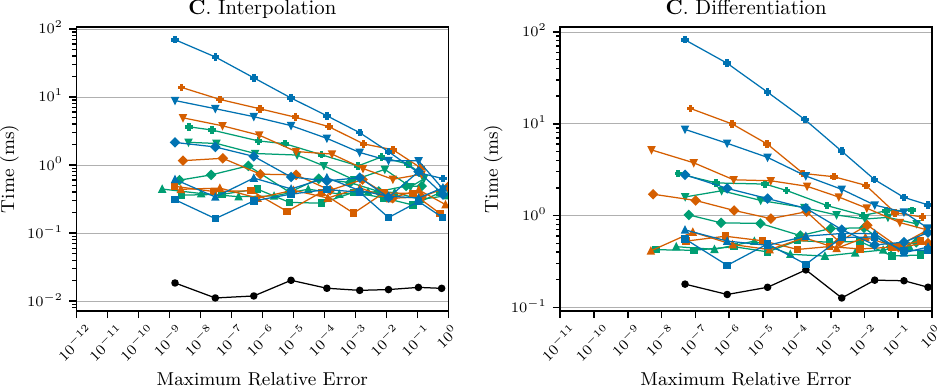}
        \caption{Anisotropic exponential function.}
        \label{fig:results_exponential}
    \end{subfigure}

    \caption[]{Results for the functions \textbf{A}, \textbf{B}, and \textbf{C}.}
    \label{fig:results_entire_functions}
\end{figure}

The numerical results for these test functions are shown in Figure~\ref{fig:results_entire_functions} and discussed in Section~\ref{subsec:discussion_numerical_results}.

\subsection{Benchmark against Chebfun and ApproxFun}
\label{subsec:runge_functions}

We further assess the numerical performance of \emph{lpFun}~\cite{Hofmann2026lpFun} version 1.1 in dimensions \(m = 1, 2, 3\) against \emph{Chebfun} \cite{Chebfun} version 5.7.0 and in dimensions \(m = 1, 2\) against \emph{ApproxFun} \cite{ApproxFun} version 0.13.
For \emph{lpFun}, we use the Euclidean \(\ell^2\)-degree, following the emphasis of \emph{Trefethen}~\cite{Trefethen2017a, Trefethen2017b} and \emph{Hecht et al.}~\cite{Hecht2026}, whereas \emph{Chebfun} and \emph{ApproxFun} use the maximal \(\ell^\infty\)-degree. In \emph{Chebfun}, the maximal degree is chosen adaptively according to a prescribed tolerance. We vary this tolerance as
\begin{equation}
    \varepsilon_k = 10^{-15(k - 1)/14}, \quad k = 1, 2, \ldots, 15.
\end{equation}
For \emph{ApproxFun}, we disable adaptive interpolation and prescribe the degree manually. We acknowledge that runtime comparisons across Python, MATLAB, and Julia are affected by language- and implementation-specific factors. The purpose of this benchmark is therefore not to provide a language-specific performance analysis, but to demonstrate the state-of-the-art algorithmic efficiency of \emph{lpFun}.

We use the following Runge-type test functions:
\begin{description}
    \item[\textbf{D}.]
    Runge function
    \begin{equation}
        f(\bm x)
        =
        \frac{1}{1 + r \|\bm x\|_2^2},
        \qquad
        r = 25.
    \end{equation}

    \item[\textbf{E}.]
    Perturbed Runge function
    \begin{equation}
        f(\bm x)
        =
        \frac{1}{1 + \sum_{i=1}^m r_{m,i} x_i^2},
        \qquad
        r_{m,i} = \frac{25}{i}.
    \end{equation}

    \item[\textbf{F}.]
    Shifted Runge-type function
    \begin{equation}
        f(\bm x)
        =
        \frac{1}{\sum_{i=1}^m (x_i-a)^2},
        \qquad
        a = \frac{5}{4}.
    \end{equation}
\end{description}

\begin{figure}[!htbp]
    \centering

    \includegraphics[width=0.9\textwidth]{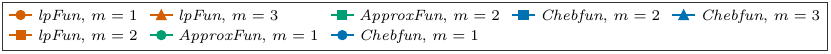}

    \vspace{0.5em}

    \begin{subfigure}{\textwidth}
        \centering
        \includegraphics[width=0.9\textwidth]{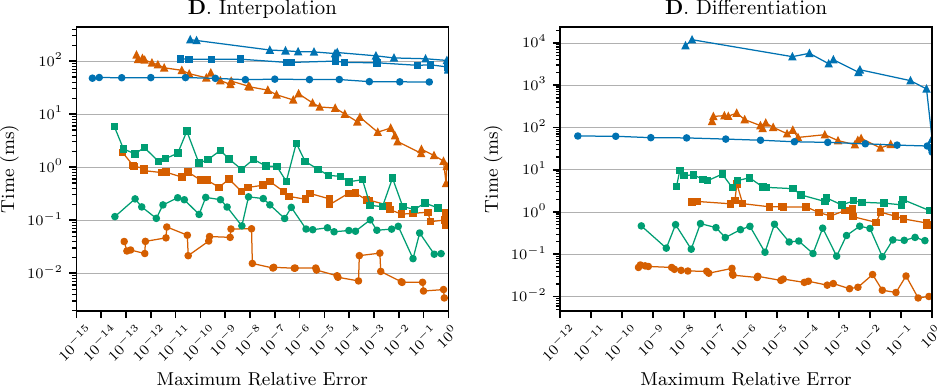}
        \caption{Runge function.}
        \label{fig:results_runge}
    \end{subfigure}

    \vspace{0.8em}

    \begin{subfigure}{\textwidth}
        \centering
        \includegraphics[width=0.9\textwidth]{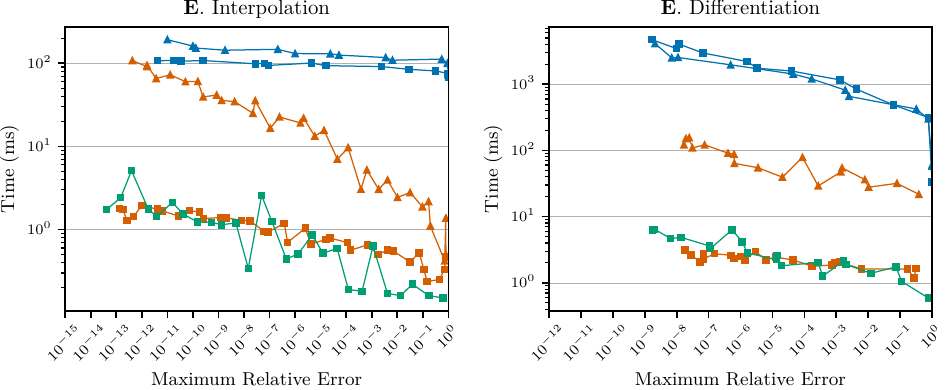}
        \caption{Perturbed Runge function.}
        \label{fig:results_perturbed_runge}
    \end{subfigure}

    \vspace{0.8em}

    \begin{subfigure}{\textwidth}
        \centering
        \includegraphics[width=0.9\textwidth]{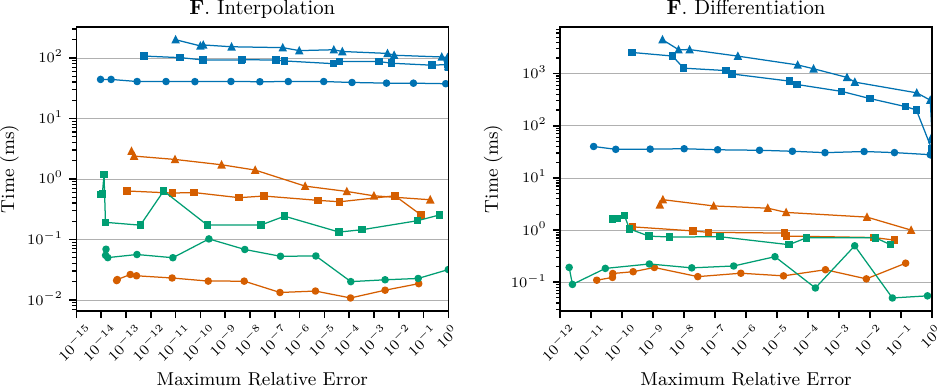}
        \caption{Shifted Runge-type function.}
        \label{fig:results_extension_runge}
    \end{subfigure}

    \caption[]{Results for the Runge-type functions \textbf{D}, \textbf{E}, and \textbf{F}.}
    \label{fig:results_runge_functions}
\end{figure}

The numerical results for these test functions are shown in Figure~\ref{fig:results_runge_functions} and discussed in Section~\ref{subsec:discussion_numerical_results}.

\pagebreak

\subsection{Discussion of the numerical results}
\label{subsec:discussion_numerical_results}
We point out the following observations.  

\paragraph{Test functions A, B, and C.} 
All three test cases exhibit the same phenomenon, that is, for dimensions \(m>1\), the runtimes are shifted by a constant margin. 
We attribute this shift primarily to the overhead introduced by enabling CPU parallelization in \emph{Numba} for \(m>1\).

For test function \textbf{A}, the \emph{Euclidean} degree is generally superior, performing best in both interpolation and differentiation in most cases. 
For test function \textbf{B}, the total and \emph{Euclidean} degrees perform similarly for interpolation, while the total degree appears to be preferable for differentiation. 
For test function \textbf{C}, the total degree dominates in interpolation and differentiation.

Overall, these results show that choosing either the total \(\ell^1\) or the \emph{Euclidean} \(\ell^2\) degree leads to significantly better performance than choosing the maximal \(\ell^\infty\) degree. 
However, even for entire functions, the preferable choice between the total and \emph{Euclidean} degree still depends on the specific function.

Although the maximal degree may occasionally yield a slightly smaller error, the corresponding accuracy gain is negligible compared with the additional computational cost.

\paragraph{Test functions D, E, and F:}

These test cases are intended as a validation of \emph{lpFun} against established packages for polynomial approximation, rather than a strict performance ranking. 

The observed results indicate that \emph{lpFun} performs competitively with \emph{ApproxFun} for the tested Runge-type examples, both for interpolation and differentiation. 
Compared with \emph{Chebfun}, \emph{lpFun} requires less runtime in these experiments, particularly as the spatial dimension increases.

In summary, the benchmarks confirm that the current implementation of the \emph{FNT} in \emph{lpFun} can match the performance of established packages in the tested regimes.

\section{Conclusions and future work}

We have constructed and analyzed the \emph{fast Newton transform} (FNT), an algorithm for \emph{Newton} interpolation in downward closed polynomial spaces \(\Pi_A\). We proved that the \emph{Newton} coefficients of the interpolant can be computed in \(\mathcal{O}(m \overline n |A|)\), where \(|A| = \dim \Pi_A\), and that evaluation of the interpolant in the underlying quasi-tensorial grid admits the same complexity. Moreover, we showed that the coefficients of the \(i\)th partial derivative of a Newton polynomial can be computed in \(\mathcal{O}(n_i |A|)\). Complementing these theoretical results, we implemented the \emph{FNT} in the open-source Python package \emph{lpFun}~\cite{Hofmann2026lpFun} for \(\ell^p\)-type index sets, making it accessible in an easy-to-use form.

Our numerical experiments demonstrate the computational performance of the \emph{FNT} algorithm and identify the total and \emph{Euclidean} degree as the more effective choices than the maximal degree, being consistent with the theoretical results of \emph{Trefethen}~\cite{Trefethen2017a}. 
Furthermore, in dimensions \(m=1,2,3\), \emph{lpFun} performs competitively with the established packages \emph{ApproxFun}~\cite{ApproxFun} and \emph{Chebfun}~\cite{Chebfun}, while requiring substantially fewer coefficients when using the \emph{Euclidean} degree.
In higher spatial dimensions \(m=4,5,6\), \emph{lpFun} has no direct counterpart among these packages and continues to exhibit the predicted computational performance.

The present work therefore provides a theoretical and practical framework for fast \emph{Newton} interpolation in downward closed polynomial spaces. 
While the \emph{Newton} basis yields a direct and canonical realization of the transform, the framework also extends to separable polynomial bases through univariate LU decompositions while retaining its time complexity of \(\mathcal O (m\overline n |A|)\). 

Future work may focus on accelerating the \emph{Chebyshev} basis case by exploiting FFT-based structure for suitable node sets that provide the nestedness required by the present framework.
Another promising direction is the extension of the present framework to low-rank and tensor-compressed representations, such as tensor train~\cite{Oseledets2011}, which could further increase the feasible dimension and polynomial degree.

\paragraph{\small Acknowledgments}
{\small
    We deeply acknowledge 
    Piotr Held,
    Damar Wicaksono,
    Leslie Greengard,
    Shidong Jiang,
    Albert Cohen,
    Oliver Sander,
    Peter Stadler,
    and Uwe Hernandez Acosta for insightful comments and helpful discussions.
    We especially thank Piotr Held for implementing the Legendre polynomial functionality in \emph{lpFun}.
}
\paragraph{\small Funding}
{\small
    This work was partially funded by the Center for Advanced Systems Understanding (CASUS), financed by Germany's Federal Ministry of Education and Research (BMBF) and by the Saxon Ministry for Science, Culture and Tourism (SMWK) with tax funds on the basis of the budget approved by the Saxon State Parliament.
}
\paragraph{\small Code availability}
{\small 
    The source code of \emph{lpFun} is available at \url{https://github.com/phil-hofmann/lpfun}.
}

\pagebreak

%Bibliography
\bibliographystyle{plain}
\bibliography{references}

\end{document}